\documentclass[11pt,a4paper]{amsart}
\usepackage[margin=3cm]{geometry}

\usepackage[T1]{fontenc}
\usepackage{amsmath,amssymb,mathtools}
\usepackage{microtype}
\usepackage[hidelinks]{hyperref}
\numberwithin{equation}{section}

\newtheorem{theorem}{Theorem}[section]
\newtheorem{proposition}[theorem]{Proposition}
\newtheorem{lemma}[theorem]{Lemma}

\theoremstyle{remark}

\newcommand{\N}{\mathcal N}
\newcommand{\ex}{\operatorname{ex}}
\newcommand{\codeg}{\operatorname{codeg}}

\newcommand{\PiFive}{\Pi_5}

\title[Pentagons in book-free graphs]
  {The Exact Number of Pentagons in Book-Free Graphs}
\author{Wentao Hu}
\address{Qiuzhen College, Tsinghua University, Beijing, China}
\email{huwt24@mails.tsinghua.edu.cn}

\begin{document}

\begin{abstract}
Let $B_k$ denote the book consisting of $k$ triangles sharing a common
edge. We prove that there exist absolute constants $\rho>0$ and $N$ such
that, for every pair of integers $n,k$ with $n\ge N$ and
$2\le k\le \rho n$, writing $n=5q+r$ with $0\le r<5$,
\[
  \ex(n,C_5,B_k)=q^{\,5-r}(q+1)^r,
\]
with the balanced complete $C_5$-blow-ups as the only extremal graphs.
The linear order of this range is best possible since for all sufficiently
large $n$, the formula fails whenever $k\ge0.14n$.
\end{abstract}

\maketitle

\section{Introduction}

For graphs \(H\) and \(G\), let \(\N(H,G)\) denote the number of
unlabelled copies of \(H\) in a graph \(G\). The \emph{generalized
Tur\'an number}
\[
  \ex(n,H,F)
  :=\max\{\N(H,G): |V(G)|=n,\ G\text{ is }F\text{-free}\}
\]
asks how many copies of \(H\) an \(n\)-vertex graph can contain without
containing \(F\) as a subgraph. Here copies are not required to be induced.
Alon and Shikhelman~\cite{AlonShikhelman2016} initiated the systematic
study of these numbers, building on earlier work on particular pairs
of graphs.

A classical instance about this framework is the problem
of maximizing the number of pentagons in a triangle-free graph. This
problem attracted great interest following a conjecture of
Erd\H{o}s that every triangle-free graph on \(n\) vertices contains at
most \((n/5)^5\) pentagons. The natural construction is a balanced
complete \(C_5\)-blow-up. More precisely, a complete \(C_5\)-blow-up
\(C_5(A_0,\ldots,A_4)\) is obtained by replacing the vertices of a
pentagon with independent sets \(A_0,\ldots,A_4\), and inserting all
edges between \(A_i\) and \(A_{i+1}\), with indices taken modulo five.
It is \emph{balanced} if its five part sizes differ by at most one.
Such a graph is triangle-free and contains exactly
\(\prod_{i=0}^4|A_i|\) pentagons. Writing \(n=5q+r\), where \(0\le r<5\),
we therefore set
\begin{equation}\label{eq:Pi}
  \PiFive(n)
  :=q^{\,5-r}(q+1)^r
  =\max_{\substack{a_0,\ldots,a_4\in\mathbb Z_{\ge0}\\
                    a_0+\cdots+a_4=n}}
     \prod_{i=0}^4 a_i.
\end{equation}

Grzesik~\cite{Grzesik2012} and, independently, Hatami, Hladk\'y,
Kr{\'a}{\v l}, Norine, and Razborov~\cite{HatamiEtAl2013} proved
Erd\H{o}s's conjectured bound
\[
  \ex(n,C_5,K_3)\le (n/5)^5.
\]
Together with the balanced blow-up construction, this determines the
asymptotic value \(n^5/3125+o(n^5)\). Hatami et al.\ also determined
the exact maximum and the extremal graphs for all sufficiently large
\(n\). Lidick\'y and Pfender~\cite[Theorem~2]{LidickyPfender2018}
subsequently completed the determination for every \(n\), proving
\[
  \ex(n,C_5,K_3)=\PiFive(n).
\]
For \(n\ge9\), the only extremal graphs are the balanced complete
\(C_5\)-blow-ups.

One can now ask how far the triangle-free condition can be relaxed
while preserving this extremal construction. Let \(B_k\) denote the
\emph{book} consisting of \(k\) triangles sharing a common edge; in
particular, \(B_1=K_3\), and \(B_2\) is the diamond. For every fixed
\(k\), the chromatic transfer bound of Gerbner and
Palmer~\cite[Theorem~2]{GerbnerPalmer2019} gives
\[
  \ex(n,C_5,B_k)\le \ex(n,C_5,K_3)+o(n^5).
\]
Combining this bound with the triangle-free construction yields
\[
  \ex(n,C_5,B_k)
  =\ex(n,C_5,K_3)+o(n^5)
  =\frac{n^5}{3125}+o(n^5).
\]

Then it is natural to ask whether the exact maximum is
also unchanged: for each fixed \(k\ge2\), is
\(\ex(n,C_5,B_k)=\PiFive(n)\) for all sufficiently large \(n\), and
are the balanced complete \(C_5\)-blow-ups the only extremal graphs?
Gerbner and Palmer explicitly listed the determination of this exact
value as an open problem in their survey~\cite[p.~45]{GerbnerPalmer2026}.
We resolve this problem and prove a stronger statement in which the
forbidden book may grow linearly with the order of the graph.

\begin{theorem}\label{thm:main}
There exist absolute constants \(\rho>0\) and \(N\) such that, for every
pair of integers \(n,k\) satisfying
\[
  n\ge N,\qquad 2\le k\le \rho n,
\]
we have
\[
  \ex(n,C_5,B_k)=\PiFive(n).
\]
Moreover, the extremal graphs are precisely the balanced complete
\(C_5\)-blow-ups.
\end{theorem}

Here the linear order of the permitted range is
best possible: Proposition~\ref{prop:linear-obstruction} shows that,
for all sufficiently large \(n\), the formula fails whenever
\(k\ge\lceil0.14n\rceil\). The requirement that \(n\) be sufficiently large cannot be dropped either.
Indeed, the triangular prism, obtained from two vertex-disjoint
triangles by adding a perfect matching between them, is diamond-free
and contains six pentagons, whereas \(\PiFive(6)=2\).

\subsection*{Proof strategy}

We first prove the following sequential formulation, and then obtain
the uniform constants in Theorem~\ref{thm:main} by a diagonal argument.

\begin{proposition}[Sequential formulation]\label{prop:sublinear}
Let \(k=k(n)\ge2\) be integer-valued and satisfy \(k=o(n)\). Then,
for all sufficiently large \(n\),
\[
  \ex(n,C_5,B_k)=\PiFive(n),
\]
and the extremal graphs are precisely the balanced complete
\(C_5\)-blow-ups.
\end{proposition}

Let \(G\) be an extremal \(B_k\)-free graph in this sublinear regime.
The proof of Proposition~\ref{prop:sublinear} proceeds from an approximate structure to an exact
one. 

First, the triangle removal lemma, together
with the triangle-free pentagon stability theorem of Pikhurko,
Slia\v{c}an, and Tyros~\cite[Theorem~1.3]{PikhurkoSliacanTyros2019},
 places \(G\) within \(o(n^2)\) edge changes of an
asymptotically balanced complete blow-up
\[
  T=C_5(A_0,\ldots,A_4).
\]

Second, for each vertex \(v\), let \(d_i(v)\) denote its degree into
\(A_i\). We derive a uniform local bound for the number of pentagons in
\(G\) containing \(v\) in terms of these five degrees.

Third, let the type class \(R_i\) consist of the vertices adjacent to
almost all vertices in both \(A_{i-1}\) and \(A_{i+1}\). We prove the following properties:
\begin{itemize}
  \item the sets \(R_0,\ldots,R_4\) are pairwise disjoint;
  \item $|R_i\mathbin{\triangle}A_i|=o(n)$, and consequently every
        consecutive pair \(R_i,R_{i+1}\) is almost completely joined;
  \item there are no edges within any \(R_i\) or between \(R_i\) and
        \(R_{i+2}\);
  \item every vertex has a type.
\end{itemize}
Consequently, \(R_0,\ldots,R_4\) partition \(V(G)\), and
\[
  G\subseteq C_5(R_0,\ldots,R_4).
\]

Finally, extremality and the integer arithmetic--geometric mean
inequality give the upper bound and force \(G\) to be a complete
\(C_5\)-blow-up whose five parts differ in size by at most one.

Section~\ref{sec:coarse} establishes the approximate structure.
Section~\ref{sec:cleaning} gives the local bounds, the construction of the type class, and provees Theorem~\ref{thm:main}.
Section~\ref{sec:obstruction} shows the
necessity of a restriction on the linear booksize range and gives a further open problem.

\section{Coarse stability}\label{sec:coarse}

We first give the coarse \(C_5\)-blow-up structure of the extremal graph.

\begin{proposition}[Coarse structure]\label{prop:coarse}
Let \(n\) tend to infinity through any unbounded sequence, let
\(k=k(n)\ge2\) satisfy \(k=o(n)\), and let \(G_n\) be an extremal
\(B_{k(n)}\)-free graph on \(n\) vertices. Then there is a complete
\(C_5\)-blow-up \(T_n=C_5(A_0,\ldots,A_4)\) on \(V(G_n)\) such that
\begin{equation}\label{eq:coarse}
  |E(G_n)\mathbin{\triangle}E(T_n)|=o(n^2),
  \qquad
  |A_i|=(1/5+o(1))n
\end{equation}
for every \(i\).
\end{proposition}

\begin{proof}
Since
\[
  3\N(K_3,G_n)\le (k(n)-1)e(G_n)=o(n^3),
\]
the triangle removal lemma~\cite[Theorem~1]{Fox2011} allows us to
delete \(o(n^2)\) edges and obtain a triangle-free graph \(F_n\).
These deletions destroy only \(o(n^5)\) pentagons, so extremality gives
\[
  \N(C_5,F_n)\ge \PiFive(n)-o(n^5).
\]
As every pentagon in a triangle-free graph is induced, the stability
theorem of Pikhurko, Slia\v{c}an, and
Tyros~\cite[Theorem~1.3]{PikhurkoSliacanTyros2019} gives a complete
\(C_5\)-blow-up \(T_n\) differing from \(F_n\), and hence from \(G_n\),
in \(o(n^2)\) edges. Consequently,
\[
  \prod_{i=0}^4|A_i|
  =\N(C_5,T_n)\ge \PiFive(n)-o(n^5).
\]
Together with \(\sum_i|A_i|=n\), asymptotic equality in the
arithmetic--geometric mean inequality yields
\(|A_i|=(1/5+o(1))n\) for every \(i\).
\end{proof}
\section{Local bounds and exact cleaning}\label{sec:cleaning}

We prove Proposition~\ref{prop:sublinear} and then
 Theorem~\ref{thm:main} in this section. Let $n$ tend to infinity
through an arbitrary unbounded sequence, let $k=k(n)\ge2$ be
integer-valued with $k=o(n)$, and for each such $n$ let $G_n$ be an
extremal $B_{k(n)}$-free graph. Proposition~\ref{prop:coarse} supplies a
complete $C_5$-blow-up
\(
  T_n=C_5(A_0,\ldots,A_4)
\)
on $V(G_n)$ satisfying \eqref{eq:coarse}. For the remainder of this
section, we suppress $n$ and simply write
\[
  G=G_n,\qquad T=T_n,\qquad k=k(n).
\]

Set
\[
  m_+:=|E(G)\setminus E(T)|,\qquad
  m_-:=|E(T)\setminus E(G)|,\qquad
  Q:=\max_i|A_i|.
\]
For every $v\in V(G)$, define
\[
  d_i(v):=|N_G(v)\cap A_i|,
\]
and
\[
  S_1(v):=\sum_i d_i(v)d_{i+1}(v),\qquad
  S_2(v):=\sum_i d_i(v)d_{i+2}(v),
\]
where all indices are taken modulo five. Proposition~\ref{prop:coarse}
and the $B_k$-free condition give the standing estimates
\begin{equation}\label{eq:standing}
  m_++m_-=o(n^2),\qquad
  Q=\left(\frac15+o(1)\right)n.
\end{equation}

For $v\in V(G)$, let $c_5(v,G)$ denote the number of pentagons of $G$
containing $v$. The following lemma bounds this quantity in terms of the
neighbourhood profile of $v$ and the template errors.

\begin{lemma}[Local path count]\label{lem:local}
For every \(v\in V(G)\),
\begin{align}
  S_1(v)
  &\le m_-+\frac{kn}2,
  \label{eq:adjacent-product}\\
  c_5(v,G)
  &\le Q^2\bigl(S_2(v)+3S_1(v)\bigr)+6m_+n^2.
  \label{eq:local-bound}
\end{align}
\end{lemma}

\begin{proof}
Every pair counted by $d_i d_{i+1}$ is either an edge of $G[N_G(v)]$
or an edge of $E(T)\setminus E(G)$. Hence
\[
  S_1(v)\le e(G[N_G(v)])+m_-.
\]
Moreover,
\[
  2e(G[N_G(v)])
  =\sum_{u\in N_G(v)}\codeg_G(u,v)
  \le kn,
\]
which proves \eqref{eq:adjacent-product}.

Deleting $v$ from a pentagon containing it leaves an unoriented
length-three path
\[
  P=x_0x_1x_2x_3
\]
in $G-v$ whose endpoints $x_0,x_3$ lie in $N_G(v)$. We distinguish two
cases.

\smallskip
\noindent\emph{Case 1: all three edges of $P$ belong to $E(T)$.}
The part indices of $x_0,x_1,x_2,x_3$ form a length-three walk in the
quotient $C_5$. Up to reversing the walk, the possibilities are as
follows, with indices taken modulo five:
\[
\begin{array}{c|c}
\text{endpoint parts} & \text{index walks}\\ \hline
A_i,\ A_{i+2}
  & (i,i-1,i-2,i-3)\\[2pt]
A_i,\ A_{i+1}
  & (i,i+1,i,i+1),\
    (i,i+1,i+2,i+1),\
    (i,i-1,i,i+1)\\[2pt]
A_i,\ A_i
  & \text{none}.
\end{array}
\]
Once the endpoints and the index walk are fixed, the two internal
vertices have at most $Q^2$ choices. Hence this case contributes at most
\[
  Q^2\bigl(S_2(v)+3S_1(v)\bigr).
\]

\smallskip
\noindent\emph{Case 2: some edge of $P$ belongs to
$E(G)\setminus E(T)$.}
Charge $P$ to one such edge. There are $m_+$ choices for the charged
edge, three possible positions for it in the path, two orientations, and
fewer than $n^2$ choices for the remaining two vertices. Therefore this
case contributes at most
\[
  6m_+n^2.
\]
Adding the two cases proves \eqref{eq:local-bound}.
\end{proof}

Before proceeding to the proof of the main theorem, we establish
one further analytic lemma. For \(x=(x_0,\ldots,x_4)\in[0,1]^5\), define

\[
  A(x)=\sum_i x_ix_{i+1},
  \qquad
  D(x)=\sum_i x_ix_{i+2}.
\]
Let
\[
  \mathcal E
  =\{e_{j-1}+e_{j+1}:j\in\mathbb Z_5\}.
\]

\begin{lemma}[Five-variable stability]\label{lem:fivevar}
\leavevmode
\begin{enumerate}
\item For every $\varepsilon>0$ there exists $\delta>0$ such that
      $A(x)\le\delta$ implies $D(x)\le1+\varepsilon$.
\item For every $\eta\in(0,1/4)$ there exist $\delta,\gamma>0$ such that,
      if $A(x)\le\delta$ and there is no $j$ satisfying
      \[
        \min\{x_{j-1},x_{j+1}\}\ge1-\eta,
      \]
      then $D(x)\le1-\gamma$.
\end{enumerate}
\end{lemma}

\begin{proof}
If $A(x)=0$, then the positive support of $x$ is an independent set in
$C_5$ and therefore has size at most two. It follows that $D(x)\le1$,
with equality precisely when $x\in\mathcal E$.

A contrary sequence for the first assertion would have a convergent
subsequence whose limit has $A=0$ and $D>1$. For the second assertion, a
contrary sequence may be chosen with $A(x^{(m)})\le1/m$ and
$D(x^{(m)})>1-1/m$. Its limit lies in $\mathcal E$, so the corresponding
two coordinates of $x^{(m)}$ are at least $1-\eta$ for all sufficiently
large $m$, a contradiction.
\end{proof}

\begin{proof}[Proof of Proposition~\ref{prop:sublinear}]

\smallskip
\noindent\emph{Step 1: Uniform local estimates and vertex types.}
For $v\in V(G)$, define
\[
  x_i(v):=\frac{d_i(v)}{Q},
  \qquad
  x(v):=(x_0(v),\ldots,x_4(v)).
\]
Thus
\[
  A(x(v))=\frac{S_1(v)}{Q^2},
  \qquad
  D(x(v))=\frac{S_2(v)}{Q^2}.
\]
Set
\[
  \alpha_n
  :=\frac{m_-}{Q^2}
    +\frac{kn}{2Q^2},
  \qquad
  r_n
  :=3\alpha_n+\frac{6m_+n^2}{Q^4}.
\]
By \eqref{eq:standing} and \eqref{eq:local-bound}, both $\alpha_n$ and $r_n$ tend to zero, and uniformly in $v$,
\begin{equation}\label{eq:A-small}
  A(x(v))\le\alpha_n,\quad   \frac{c_5(v,G)}{Q^4}
  \le D(x(v))+r_n.
\end{equation}

Fix $\eta\in(0,1/4)$. For each $j\in\mathbb Z_5$, define
\begin{equation}\label{eq:type}
  R_j
  :=\bigl\{v\in V(G):
      \min\{x_{j-1}(v),x_{j+1}(v)\}\ge1-\eta
    \bigr\}.
\end{equation}
We call the vertices in $R_j$ vertices of \emph{type $j$}, and set
\[
  U:=V(G)\setminus\bigcup_{j\in\mathbb Z_5}R_j.
\]
The vertices in $U$ are called \emph{untyped}.

For all sufficiently large $n$, the sets $R_0,\ldots,R_4$ are pairwise
disjoint. Indeed, if $v\in R_i\cap R_j$ for distinct $i,j$, then two
adjacent coordinates of $x(v)$ are both at least $1-\eta$. Hence
\[
  A(x(v))\ge(1-\eta)^2,
\]
contrary to \eqref{eq:A-small}.

Lemma~\ref{lem:fivevar} and \eqref{eq:A-small} imply that there
exist a constant $\gamma>0$ and a sequence $\varepsilon_n\to0$ such that
\begin{align}
  c_5(v,G)
  &\le(1+\varepsilon_n)Q^4
  &&\text{for every }v\in V(G),
  \label{eq:all-local}\\
  c_5(v,G)
  &\le(1-\gamma)Q^4
  &&\text{for every }v\in U.
  \label{eq:untyped-gap}
\end{align}
Here $\gamma$ has been decreased, if necessary, to absorb the error
$r_n$.

\smallskip
\noindent\emph{Step 2: Matching the type classes to the coarse partition.}
Define the error degree of $v$ by
\[
\operatorname{err}(v)
:=|N_G(v)\mathbin{\triangle}N_T(v)|.
\]
Put
\[
  M_n:=m_++m_-,
  \qquad
  \beta_n:=\frac{\sqrt{M_n}+1}{n}=o(1).
\]
If $v\in A_i$ and $\operatorname{err}(v)\le\beta_n n$, then
\[
  d_{i\pm1}(v)
  \ge |A_{i\pm1}|-\beta_n n
  =(1-o(1))Q
  >(1-\eta)Q,
\]
and hence $v\in R_i$. Conversely, let $v\in R_i\setminus A_i$, and let
$j\ne i$ be the index for which $v\in A_j$. If
$\operatorname{err}(v)\le\beta_n n$, then the same estimate gives
\[
  d_{j\pm1}(v)
  \ge |A_{j\pm1}|-\beta_n n
  >(1-\eta)Q,
\]
so $v\in R_j$, contradicting the pairwise disjointness of
$R_0,\ldots,R_4$. Consequently,
\[
  R_i\mathbin{\triangle}A_i
  \subseteq
  \{v:\operatorname{err}(v)>\beta_n n\}.
\]
Since
\[
  \sum_{v\in V(G)}\operatorname{err}(v)=2M_n,
\]
Markov's inequality and \eqref{eq:standing} give
\[
  \bigl|\{v:\operatorname{err}(v)>\beta_n n\}\bigr|
  \le \frac{2M_n}{\sqrt{M_n}+1}
  \le 2\sqrt{M_n}
  =o(n).
\]
This implies
\[
  |R_i\mathbin{\triangle}A_i|=o(n),
  \qquad
  |R_i|=\left(\frac15+o(1)\right)n.
\]

We next show that $R_i$ and $R_{i+1}$ are almost completely joined. Fix
$i$, and consider a non-edge $uv$ with $u\in R_i$ and
$v\in R_{i+1}$. If $u\in A_i$ and $v\in A_{i+1}$, then
$uv\in E(T)\setminus E(G)$, so this non-edge is counted by $m_-$. Otherwise,
\[
  u\in R_i\setminus A_i
  \qquad\text{or}\qquad
  v\in R_{i+1}\setminus A_{i+1}.
\]
Consequently,
\[
  \begin{aligned}
  |R_i||R_{i+1}|-e_G(R_i,R_{i+1})
  &\le m_-+n|R_i\setminus A_i|
          +n|R_{i+1}\setminus A_{i+1}|\\
  &=o(n^2).
  \end{aligned}
\]

\smallskip
\noindent\emph{Step 3: Embedding $G$ into a $C_5$-blow-up.}
We claim that, for every $i$,
\begin{equation}\label{eq:forbidden-type-pairs}
  E(G[R_i])=\varnothing,
  \qquad
  E_G(R_i,R_{i+2})=\varnothing.
\end{equation}
Indeed, two vertices of type $i$ have at least
\[
  2(1-\eta)Q-|A_{i-1}|
  \ge(1-2\eta)Q
\]
common neighbours in $A_{i-1}$. Similarly, vertices of types $i$ and
$i+2$ have at least $(1-2\eta)Q$ common neighbours in $A_{i+1}$. By
\eqref{eq:standing},
\[
  (1-2\eta)Q>k
\]
for all sufficiently large $n$. Hence no such pair can be an edge. This proves \eqref{eq:forbidden-type-pairs}.
In particular,
\begin{equation}\label{eq:independent-type-neighbourhood}
  R_{i-1}\cup R_{i+1}
  \quad\text{is independent for every }i.
\end{equation}

Now suppose that an untyped vertex $v\in U$ remains. Fix an index $i$,
and let $G'$ be obtained from $G$ by deleting $v$ and adding a new vertex
$z$ with
\[
  N_{G'}(z)=R_{i-1}\cup R_{i+1}.
\]
By \eqref{eq:independent-type-neighbourhood}, the neighbourhood of $z$ is
independent. Thus $z$ lies in no triangle, and no old edge gains $z$ as
a common neighbour. It follows that $G'$ is still $B_k$-free.

Consider the four-part paths
\[
  R_{i-1}-R_{i-2}-R_{i+2}-R_{i+1}.
\]
By Step~2, the total number $\ell_i$ of missing edges in the three
displayed consecutive pairs is $o(n^2)$. Every missing edge excludes at
most $n^2$ choices of a four-vertex path. Therefore
\begin{align*}
  c_5(z,G')
  &\ge \prod_{j\ne i}|R_j|-\ell_i n^2\\
  &=(1-o(1))Q^4.
\end{align*}
For all sufficiently large $n$, this is strictly larger than
$(1-\gamma)Q^4$, whereas \eqref{eq:untyped-gap} gives
\[
  c_5(v,G)\le(1-\gamma)Q^4.
\]
Since every pentagon not containing the replaced vertex is unchanged,
\[
  \N(C_5,G')-\N(C_5,G)
  =c_5(z,G')-c_5(v,G)>0,
\]
contradicting the extremality of $G$. Hence $U=\varnothing$, so
$R_0,\ldots,R_4$ is a partition of $V(G)$. Together with
\eqref{eq:forbidden-type-pairs}, this gives
\[
  G\subseteq T':=C_5(R_0,\ldots,R_4).
\]

\smallskip
\noindent\emph{Step 4: Completion and balancing.}
The balanced complete $C_5$-blow-up is $B_k$-free. Hence, by the
extremality of $G$ and the integer form of the arithmetic--geometric mean
inequality,
\[
  \PiFive(n)
  \le \N(C_5,G)
  \le \N(C_5,T')
  =\prod_{i=0}^4 |R_i|
  \le \PiFive(n).
\]
Thus equality holds throughout. Equality in the final inequality implies
that the five part sizes differ by at most one, while
$\N(C_5,G)=\N(C_5,T')$ forces $G=T'$: all five parts are nonempty for
sufficiently large $n$, and every edge of $T'$ lies in a transversal copy
of $C_5$. This proves Proposition~\ref{prop:sublinear}.
\end{proof}

\begin{proof}[Proof of Theorem~\ref{thm:main}]
Suppose that no such constants \(\rho,N\) exist. Recursively choose
integers
\[
  n_j>\max\{n_{j-1},2j\},
  \qquad
  2\le k_j\le n_j/j,
\]
for which either the asserted extremal value or the asserted
classification of extremal graphs fails, where \(n_0:=0\).
Define an integer-valued function on all positive integers by
\[
  \widetilde k(n):=
  \begin{cases}
    k_j,&n=n_j\text{ for some }j,\\
    2,&\text{otherwise}.
  \end{cases}
\]
Then \(\widetilde k(n)/n\longrightarrow0\), but the conclusion of
Proposition~\ref{prop:sublinear} fails at every \(n_j\).
This contradiction proves the theorem.
\end{proof}

\section{A linear obstruction}\label{sec:obstruction}

The following construction shows that the conclusion of
Theorem~\ref{thm:main} cannot extend to booksize density $0.14$.

\begin{proposition}\label{prop:linear-obstruction}
For every fixed $c\ge 0.14$ and all sufficiently large $n$, there exists
a $B_{\lfloor cn\rfloor}$-free $n$-vertex graph $H$ such that
\[
  \N(C_5,H)>\PiFive(n).
\]
Consequently, balanced complete $C_5$-blow-ups are not extremal in this
range.
\end{proposition}

\begin{proof}
Set
\[
  c_0:=0.14,\qquad
  s:=\lfloor c_0n\rfloor-1,\qquad
  t:=\left\lfloor\frac{n-6s}{3}\right\rfloor,
\]
and let
\[
  H:=2K_{s,s,s}\sqcup K_{t,t,t}
     \sqcup (n-6s-3t)K_1.
\]
 The graph $H$ has booksize $s$ for all sufficiently large $n$.
Thus $H$ is 
$B_{\lfloor cn\rfloor}$-free.

Moreover,
\[
  \N(C_5,H)
  =6s^3(s-1)^2+3t^3(t-1)^2>\left(
      \frac{1}{3125}+\frac{1}{10^6}
    \right)n^5.
\]

we obtain $\N(C_5,H)>\PiFive(n)$ for all sufficiently large $n$.
\end{proof}

\smallskip
\noindent\emph{Open problem.}
More generally, it would be interesting to determine the function
\[
  F(c):=
  \limsup_{n\to\infty}
  \frac{\ex(n,C_5,B_{\lfloor cn\rfloor})}{n^5},
  \qquad c>0.
\]
Theorem~\ref{thm:main} gives
\[
  F(c)=\frac1{3125}
\]
for every sufficiently small fixed $c>0$, whereas
Proposition~\ref{prop:linear-obstruction} shows that
\[
  F(c)>\frac1{3125}
  \qquad\text{for every }c\ge0.14.
\]
Determining the exact transition point, the function $F(c)$ beyond this
point, and the corresponding extremal constructions remain open.

\section*{AI assistance}

The main mathematical ideas of this paper are due to the author.
During the early exploratory stages of this work, the author
discussed the problem with generative AI tools from OpenAI.
These tools were also used to assist with literature searches
and proof checking.
The author takes full responsibility for the content of this paper.

\bibliographystyle{amsplain}
\bibliography{references_submit}

\end{document}